\documentclass[11pt,amsfonts]{article}
\usepackage{graphicx}
\usepackage{latexsym}
\usepackage{amssymb}
\usepackage{amsmath}
\usepackage{enumerate}

\usepackage{layout}
\usepackage{eufrak}
\newtheorem{prop}{Proposition}
\newtheorem{lemma}{Lemma}

\newtheorem{theorem}{Theorem}
\newtheorem{remark}{Remark}

\def\real{{\mathord{{\rm I\kern-2.8pt R}}}}        
\def\inte{{\mathord{{\rm I\kern-2.8pt N}}}}

\def\sZZ{{\rm Z\kern-2.8ptem{}Z}}

\def\z{{\mathchoice
  {\sZZ}
  {\sZZ}
  {\rm Z\kern-0.30em{}Z}
  {\rm Z\kern-0.25em{}Z} }}
\def\sQQ{{\kern 0.27em \vrule height1.45ex width0.03em depth0em
          \kern-0.30em \rm Q}}
\def\qu{{\mathchoice
    {\sQQ}
    {\sQQ}
  {\kern 0.225em \vrule height1.05ex width0.025em depth0em \kern-0.25em \rm Q}
  {\kern 0.180em \vrule height0.78ex width0.020em depth0em \kern-0.20em \rm Q}
        }}
\def\sCC{{\kern 0.27em \vrule height1.45ex width0.03em depth0em
          \kern-0.30em \rm C}}
\def\complex{{\mathchoice
    {\sCC}
    {\sCC}
  {\kern 0.225em \vrule height1.05ex width0.025em depth0em \kern-0.25em \rm C}
  {\kern 0.180em \vrule height0.78ex width0.020em depth0em \kern-0.20em \rm C}
        }}

\newcommand{\ba}{\begin{array}}
\newcommand{\ea}{\end{array}}
\newcommand{\be}{\begin{equation}}
\newcommand{\ee}{\end{equation}}
\newcommand{\bea}{\begin{eqnarray}}
\newcommand{\eea}{\end{eqnarray}}
\newcommand{\beaa}{\begin{eqnarray*}}
\newcommand{\eeaa}{\end{eqnarray*}}

\newcommand{\eps}{\varepsilon}

\def\z{\zeta}

\font\tenmath=msbm10 \font\sevenmath=msbm7 \font\fivemath=msbm5
\newfam\mathfam \textfont\mathfam=\tenmath
\scriptfont\mathfam=\sevenmath \scriptscriptfont\mathfam=\fivemath

\def \={{\buildrel {\rm (law)} \over =}}

\def\qed{ \hfill \vrule width.25cm height.25cm depth0cm\smallskip}

\newcommand{\basa}{\begin{assumption}}
\newcommand{\easa}{\end{assumption}}

\newcommand{\bas}{\begin{assum}}
\newcommand{\eas}{\end{assum}}

\def\liminf{\mathop{\underline{\rm lim}}}

\newcommand{\ignore}[1]{}
\begin{document}

\renewcommand{\thefootnote}{\fnsymbol{footnote}}

\renewcommand{\thefootnote}{\fnsymbol{footnote}}

\title{Least squares estimator for the parameter of the fractional Ornstein-Uhlenbeck sheet}
\author{Jorge Clarke De la Cerda $^{1,}$ \footnote{Partially supported by the MECESUP proyect UCO-0713 and the CONICYT-ECOS program C10E03.} $\qquad $
Ciprian A. Tudor $^{2,}$ \footnote{Associate member of the team Samm, Universit\'e de Paris 1 Panth\'eon-Sorbonne. Partially supported by the ANR grant "Masterie" BLAN 012103. }\vspace*{0.1in} \\
$^{1}$ CI$^{2}$MA, Departamento de Ingenier\'ia Matem\'atica, Universidad de Concepci\'on,\\
Casilla 160-C, Concepci\'on, Chile. \\
jclarke@udec.cl \vspace*{0.1in} \\
 $^{2}$ Laboratoire Paul Painlev\'e, Universit\'e de Lille 1\\
 F-59655 Villeneuve d'Ascq, France.\\
 \quad tudor@math.univ-lille1.fr\vspace*{0.1in}}

\maketitle

\begin{abstract}
We will study the least square estimator $\widehat{\theta }_{T,S}$ for the drift parameter $\theta$ of the fractional Ornstein-Uhlenbeck sheet which is defined as the solution of the Langevin  equation
\begin{equation*}
    X_{t,s}= -\theta \int^{t}_{0} \int^{s}_{0} X_{v,u}dvdu + B^{\alpha, \beta}_{t,s}, \qquad (t,s) \in [0,T]\times [0,S].
\end{equation*}
 driven by the fractional Brownian sheet $B^{\alpha ,\beta}$  with Hurst parameters $\alpha, \beta$ in $(\frac{1}{2}, \frac{5}{8})$. Using the properties of multiple Wiener-It\^o integrals we prove that the estimator is strongly consistent for the parameter $\theta$. In contrast to the one-dimensional case, the estimator $\widehat{\theta}_{T,S}$ is not asymptotically normal.
\end{abstract}

\vskip0.3cm

{\bf 2010 AMS Classification Numbers:}  60G15, G0H07, 60G35, 62M40.

 \vskip0.3cm


\section{Introduction}
We consider the two-parameter fractional Ornstein-Uhlenbeck process defined as the solution of the stochastic equation

\begin{equation}\label{1}
    X_{t,s}= -\theta \int^{t}_{0} \int^{s}_{0} X_{v,u}dvdu + B^{\alpha, \beta}_{t,s}, \qquad (t,s) \in [0,T]\times [0,S].
\end{equation}
Here $B^{\alpha , \beta}$ denotes a fractional Brownian sheet with Hurst parameters $\alpha , \beta \in (\frac{1}{2}, 1)$. We also suppose that $X_{0,0}=X_{t,0}=X_{0,s}=0$ for every $t,s$.  Our goal is to estimate the unknown parameter $\theta$ from the continuous time observation of the solution $(X_{t,s}) _{(t,s)\in [0,T]\times [0,S]}$.

The development of the stochastic analysis for
fractional Brownian motion (fBm) naturally led to the study of the statistical inference for stochastic equations driven by this process. There already exists an important literature related to these aspects. We refer, among others to \cite{BTT1}, \cite{BTT2}, \cite{HN}, \cite{KB},  \cite{TV}. Statistical analysis  of the stochastic differential equations (SDE) driven by the fractional Brownian sheet has been less considered. We refer to the paper \cite{SoTu} for the study of the maximum likelihood estimator for a SDE with additive fractional Brownian sheet noise (see also \cite{DK} or \cite{APZ} for the case when the noise is a standard Brownian sheet).

In this paper we propose a least square estimator for the unknown parameter $\theta$ following the approach in \cite{HN}. This estimator is obtained by formally minimizing with respect to $\theta$ the expression
\begin{equation*}
\int_{0}^{T} \int_{0}^{S} \left| \frac{\partial ^{2}}{\partial t
\partial s } X_{t,s} + \theta X_{t,s} \right| ^{2} dsdt.
\end{equation*}
We obtain the following estimator
\begin{equation}\label{est1}
    \widehat{\theta}_{T,S}=  - \frac{\int _{0}^{T}
\int_{0}^{S}X_{t,s}dX_{t,s}} {\int _{0}^{T}
\int_{0}^{S}X^{2}_{t,s}dtds}.
\end{equation}
The integral with respect to $dX_{t,s} $ is understood as the sum of the standard Lebesgue integral $-\theta\int_{0}^{T} \int_{0}^{S} dtds X_{t,s} ^{2}$ and of the stochastic integral $\int_{0} ^{T} \int_{0} ^{S} X_{t,s} dB_{t,s}^{\alpha , \beta}$ which is a divergence type integral with respect to the fractional Brownian sheet $B^{\alpha , \beta}$ (it will be defined in Section 2, we also refer to \cite{K}, \cite{K1}, \cite{K2}, \cite{TV2}, \cite{TV3} for the stochastic integration with respect to $B^{\alpha, \beta}$).
Using (\ref{1}) and (\ref{est1}) we can write
\begin{equation}\label{theta-theta}
\widehat{\theta}_{T,S}-\theta =- \frac{\int _{0}^{T}
\int_{0}^{S}X_{t,s}dB^{\alpha, \beta}_{t,s}} {\int _{0}^{T}
\int_{0}^{S}X^{2}_{t,s}dtds}.
\end{equation}
We will study the asymptotic behavior of the least square estimator $\widehat{\theta}_{T,S}$ as $T,S \to \infty$. Our tools are the multiple stochastic integrals and the Malliavin calculus. Actually, the nominator and the denominator of the right hand side of (\ref{theta-theta}) can be expressed as multiple integrals of order 2 with respect to the fractional Brownian sheet and from this we will obtain concrete estimates for their moments. We will prove that the estimator (\ref{est1}) is a strongly consistent estimator in the sense that it converges almost surely to the true value of the parameter $\theta$. This result is  similar to the one-dimensional case (see \cite{HN}), however the approach presented in \cite{HN} is not possible to be followed for the two-parameter case, instead we use among other tools, the hypercontractivity of multiples integrals. By contrary, in the two-parameter case, the least square estimator does not preserve the asymptotic normality as in  the one-parameter case. This will be noticed at the end of our work by using criteria  for the asymptotic normality of sequences of multiple integrals in terms of Malliavin calculus.

Our paper is structured as follows. Section 2 contains some preliminaries on multiple integrals and fractional Brownian sheet. In Section 3 we discuss the relation between the solution to (\ref{1}) and the Bessel function of order 0. Section 4 contains the proof of the consistency of the least square estimator while Section 5 is devoted to a discussion about the asymptotic normality of the estimator.

\section{Preliminaries}
\setcounter{theorem}{0}

Let us introduce the elements from stochastic analysis that we will
need in the paper. Consider ${\mathcal{H}}$ a real separable Hilbert
space and $(B (\varphi), \varphi\in{\mathcal{H}})$ an isonormal
Gaussian process on a probability space $(\Omega, {\cal{A}}, P)$, that is, a centered Gaussian family of random
variables such that $\mathbf{E}\left( B(\varphi) B(\psi) \right)  =
\langle\varphi, \psi\rangle_{{\mathcal{H}}}$.
Denote by $I_{n}$ the multiple stochastic integral with respect to
$B$ (see \cite{N}). This $I_{n}$ is actually an isometry between the
Hilbert space ${\mathcal{H}}^{\odot n}$(symmetric tensor product)
equipped with the scaled norm
${\sqrt{n!}}\Vert\cdot\Vert_{{\mathcal{H}}^{\otimes n}}$ and
the Wiener chaos of order $n$ which is defined as the closed linear
span of the random variables $H_{n}(B(\varphi))$ where
$\varphi\in{\mathcal{H}}, \Vert\varphi\Vert_{{\mathcal{H}}}=1$ and
$H_{n}$ is the Hermite polynomial of degree $n\geq 1$
\begin{equation*}
H_{n}(x)=\frac{(-1)^{n}}{n!} \exp \left( \frac{x^{2}}{2} \right)
\frac{d^{n}}{dx^{n}}\left( \exp \left( -\frac{x^{2}}{2}\right)
\right), \hskip0.5cm x\in \mathbb{R}.
\end{equation*}

 The isometry of multiple integrals can be written as: for $m,n$ positive integers,
\begin{eqnarray}
\mathbf{E}\left(I_{n}(f) I_{m}(g) \right) &=& n! \langle f,g\rangle _{{\mathcal{H}}^{\otimes n}}\quad \mbox{if } m=n,\nonumber \\
\mathbf{E}\left(I_{n}(f) I_{m}(g) \right) &= & 0\quad \mbox{if } m\not=n. \label{iso}
\end{eqnarray}
It also holds that
\begin{equation*}
I_{n}(f) = I_{n}\big( \tilde{f}\big)
\end{equation*}
where $\tilde{f} $ denotes the symmetrization of $f$ defined by $$\tilde{f}%
(x_{1}, \ldots , x_{x}) =\frac{1}{n!} \sum_{\sigma \in {\cal S}_{n}}
f(x_{\sigma (1) }, \ldots , x_{\sigma (n) } ). $$
The Malliavin derivative acts on multiple integrals $F=I_{n}(f)$ in the following way: for every $s$
\begin{equation*}
D_{s} I_{n} = nI_{n-1} (f(\cdot , s))
\end{equation*}
where $" \cdot "$ above denotes $n-1$ variables. We recall the following hypercontractivity property  for the $L^{p}$ norm of a multiple
stochastic integral (see \cite[Theorem 4.1]{Major})
\begin{equation}\label{hyper}
  \mathbf{E} \left| I_{m}(f) \right| ^{2m} \leq c_{m} \left( \mathbf{E} I_{m}(f)^{2}
  \right)^{m}
\end{equation}
where $c_{m}$ is an explicit positive constant and $f\in{\cal{H}}^{\otimes m}$.

In this work we use Malliavin calculus and multiple integrals with respect to the fractional Brownian sheet (fBs). Let us define this process and its associated Hilbert space. The fBs with Hurst parameters $ \alpha, \beta \in (0,1), (B^{\alpha,\beta}_{t,s}, t,s \in [0,T] \times [0,S]) $ is a zero mean Gaussian process with covariance

\begin{eqnarray}
    \mathbf{E}\left( B^{\alpha,\beta}_{t,s}, B^{\alpha,\beta}_{u,v} \right) &=& \mathcal{R}^{\alpha}(t,u)\mathcal{R}^{\beta}(s,v) \nonumber \\
    &:=& \frac{1}{2}\left( t^{2\alpha}+u^{2\alpha}-\vert t-u \vert ^{2\alpha} \right) \frac{1}{2}\left( s^{2\beta}+v^{2\beta}-\vert s-v \vert ^{2\beta} \right) \label{cov}
\end{eqnarray}
given for all $t,u \in [0,T]^{2}$ and $s,v \in [0,S]^{2}$.\\

We assume that $B^{\alpha,\beta}$ is defined on a complete probability space $\left( \Omega, \mathcal{A}, \mathbb{P} \right) $ such that $\mathcal{A}$ is generated by $B^{\alpha,\beta}$. Fix a time interval $[0,T]\times [0,S]$, denote by $\xi$ the set of real valued step functions on $[0,T]\times [0,S]$ and let $\mathcal{H}^{\alpha, \beta}$ be the Hilbert space defined as the closure of $\xi$ with respect to the scalar product

\begin{equation*}
 \langle 1_{[0,t]\times[0,s]} , 1_{[0,u]\times[0,v]} \rangle _{\mathcal{H}^{\alpha, \beta}} =
  \mathcal{R}^{\alpha}(t,u)\mathcal{R}^{\beta}(s,v)
\end{equation*}
where $\mathcal{R}^{\alpha}(t,u)\mathcal{R}^{\beta}(s,v)$ is the covariance function of the fBs, given in (\ref{cov}). The mapping $1_{[0,t]\times[0,s]} \longmapsto B^{\alpha,\beta}_{t,s}$ can be extended to a linear isometry between $\mathcal{H}^{\alpha, \beta}$ and the Gaussian space $\mathcal{H}^{\alpha, \beta}_{1}$ spanned by $B^{\alpha,\beta}$ which is a closed subspace of $L^{2}(\Omega,\mathcal{A},\mathbb{P})$. We denote this isometry by $\varphi \longmapsto B^{\alpha,\beta}(\varphi)$.
Fix $\alpha, \beta > \frac{1}{2}$, in this case we have that for every $f,g \in \mathcal{H}^{\alpha, \beta}$ the scalar product has the form
\begin{equation} \label{pe}
\left\langle f,g \right\rangle _{\mathcal{H}^{\alpha, \beta}} = c(\alpha)c(\beta) \int^{T}_{0} \int^{S}_{0} \int^{T}_{0} \int^{S}_{0} f(a,b)g(m,n) \vert a-m \vert ^{2\alpha -2} \vert b-n \vert ^{2\beta -2} dadbdmdn
\end{equation}
and $c(\alpha)=\alpha (2\alpha -1)$.


\section{About the solution}

The equation  (\ref{1}) has been studied in several papers (see \cite{ENO}, \cite{NT}). It has been showed that  for $\theta >0$ and $\alpha, \beta >\frac{1}{2}$ equation (\ref{1})  admits an unique strong solution which can be expressed as
\begin{equation}\label{xf}
    X_{t,s}= \int^{T}_{0} \int_{0}^{S} f(t,s,t_{0},s_{0}) dB^{\alpha, \beta}_{t_{0},s_{0}}
\end{equation}
where
\begin{equation}\label{f}
    f(t,s,t_{0},s_{0})= 1_{[0,t]}(t_{0}) 1_{[0,s]}(s_{0}) \sum_{n\geq 0}(-1)^{n} \theta^{n} \frac{(t-t_{0})^{n} (s-s_{0})^{n}}{(n!)^{2}}.
\end{equation}
We will call the process $X$ solution to (\ref{1}) as the fractional Ornstein-Uhlenbeck sheet. It is a Gaussian process since it is given by a multiple integral of order 1 (Wiener integral actually) with respect to the Gaussian process $B^{\alpha, \beta}$.  We mention that the solution to (\ref{1}) behaves differently as its one-dimensional counterpart which is fractional Ornstein-Uhlenbeck process introduced in \cite{CKM}, this will make our analysis quite different from \cite{HN}. For example, we note that the solution of some stochastic differential equations driven by the Brownian sheet or fractional Brownian sheet, which are positive in the one-parameter case,  can take negative values with strictly positive probability (see \cite{NT} or \cite{N}).

A key element of our analysis is the fact that the solution $X$ (more precisely the kernel $f$ of the solution) can be expressed is terms of the Bessel function of the first kind.  Let us consider the Bessel function of order 0 given, for every $x\in \mathbb{R}$, by
\begin{equation*}
J_{0}(x)=\sum_{n\geq 0} \frac{(-1) ^{n}}{n! ^{2}} \left( \frac{x}{2} \right) ^{2n}
\end{equation*}
This Bessel function admits the integral representation, for every $x\in \mathbb{R}$
\begin{equation*}
J_{0} (x)= \frac{1}{\pi } \int_{0} ^{\pi } \cos \left( x\sin \rho \right) d\rho.
\end{equation*}
The kernel $f$ in  (\ref{f}) of the solution $(X_{t,s}) _{t,s\in [0,T] \times [0,S]}$ can be expressed as
\begin{eqnarray}
f(t,s, u, v)&=& 1_{[0,t]}(u) 1_{[0,s] }(v) J_{0} \left( 2 \sqrt
{\theta (t-u) (s-v)} \right) \label{fjo} \\
&=& 1_{[0,t]}(u) 1_{[0,s] }(v) \frac{1}{\pi } \int_{0} ^{\pi } \cos \left(2\sqrt {\theta (t-u)
(s-v)} \sin \rho \right) d\rho.\nonumber
\end{eqnarray}
Let us also recall the following property of the Bessel function (see e.g. \cite{E}) which will play an important role for our estimates: for $x$ large enough
\begin{equation}\label{asymbessel}
J_{0} (x)\sim \sqrt{ \frac{2}{\pi x}} \cos (x-\frac{\pi}{4})
\end{equation}
(the symbol $\sim$ means that the two sides have the same limit as $x\to \infty$).

\section{Asymptotic behavior of the least square estimator}
\setcounter{theorem}{0}

In this section we study the asymptotic behavior of the estimator $\widehat{\theta} _{T,S}$ defined in (\ref{est1}). More precisely, we will show that this estimator is strongly consistent for the parameter $\theta$, that is, $\widehat{\theta}_{T,S}$ converges to $\theta $ almost surely as $T,S \to \infty$. To this end we will analyze separately the nominator and the denominator appearing in the right hand side of the expression (\ref{theta-theta}). Let us start with the study of the nominator. It can be written as the double stochastic integral
\begin{equation}\label{FTS}
 \int_{0} ^{T} \int_{0}^{S} X_{t,s} dB^{\alpha , \beta}_{t,s} :=   F_{T,S}:=I_{2}\left( f(u,v,t,s) \right)
\end{equation}
where the kernel $f$ is given by (\ref{f}) and the integral $I_{2}$ acts with respect to the variables $(u,v), (t,s)$.

We will estimate first the $L^{2}$ norm of $F_{T,S}$. We have the following result.

\begin{prop}\label{pn1}
For every $\eps >0$ and for $\alpha , \beta \in (\frac{1}{2}, \frac{5}{8})$,
\begin{equation}\label{z1}
\mathbf{E}\left( T ^{-2\alpha + \frac{1}{4}-\eps } S^{-2\beta + \frac{1}{4}-\eps} \int_{0}^{T} \int_{0}^{S} X_{t,s} dB^{\alpha, \beta}_{t,s} \right) ^{2} \to 0 \mbox{ when } T,S \to \infty.
\end{equation}
Moreover for $T,S$ large enough we have
$$\mathbf{E}\left( T ^{-2\alpha + \frac{1}{4} } S^{-2\beta + \frac{1}{4}} \int_{0}^{T} \int_{0}^{S} X_{t,s} dB^{\alpha, \beta}_{t,s} \right) ^{2}<C$$
where $C$ is a strictly positive constant not depending on $T,S$.
\end{prop}
{\bf Proof: }We calculate the $L^{2}$ norm of the random variable $I_{2}\left( f(u,v,t,s) \right)$. By the isometry property of multiple integrals (\ref{iso}) and since $\Vert \tilde{f}\Vert _{({\cal{H}}^{\alpha ,\beta })^{\otimes 2 }}\leq \Vert f\Vert _{({\cal{H}}^{\alpha ,\beta })^{\otimes 2 }}$ this norm can be handles as follows
\begin{eqnarray*}
    I_{T,S} &\leq & \int_{[0,T]^{4}}dtdt_{0} dudu_{0} \int_{[0,S] ^{4}}dsds_{0} dv dv_{0}\\
&& \times f(t,s,u,v)f(t_{0},s_{0},u_{0},v_{0})\vert u-u_{0} \vert^{2\alpha-2} \vert v-v_{0} \vert^{2\beta-2} \vert t-t_{0} \vert^{2\alpha-2} \vert s-s_{0} \vert^{2\beta-2} \\
&=&   \int_{0} ^{T} dt \int_{0}^{t} du \int_{0} ^{T} dt_{0} \int_{0} ^{t_{0}} du_{0}\int_{0} ^{S} ds \int_{0}^{s} dv \int_{0} ^{S} ds_{0} \int_{0} ^{s_{0} } dv_{0} \\
&&\times J_{0} \left( 2\sqrt{ \theta (t-u)(s-v)} \right) J_{0} \left(2 \sqrt{\theta (t_{0} -u_{0}) (s_{0} -v_{0} ) } \right) \\
&&\times \vert t-t_{0} \vert ^{2\alpha -2} \vert u-u_{0} \vert ^{2\alpha -2} \vert s-s_{0} \vert ^{2\beta -2} \vert v-v_{0} \vert ^{2\beta -2}.
\end{eqnarray*}
By making the change of variables $\tilde{t}= \frac{t}{T}$, $\tilde{u}= \frac{u}{T}$ and similarly for the other variables, we obtain
\begin{eqnarray}
I_{T,S}&=& T^{4\alpha -4} T^{4} S^{4\beta -4} S^{4} \int_{0} ^{1}dt \int_{0}^{t} du\int_{0} ^{1}dt_{0} \int_{0}^{t_{0}} du_{0} \int_{0} ^{1}ds \int_{0}^{s} dv \int_{0}^{1} ds_{0} \int_{0}^{s_{0}} dv_{0} \nonumber\\
&&\times J_{0} \left( 2\sqrt{ \theta (t-u)(s-v)TS} \right) J_{0} \left(2 \sqrt{\theta (t_{0} -u_{0}) (s_{0} -v_{0} )TS } \right) \nonumber \\
 &&\times \vert t-t_{0} \vert ^{2\alpha -2} \vert u-u_{0} \vert ^{2\alpha -2} \vert s-s_{0} \vert ^{2\beta -2} \vert v-v_{0} \vert ^{2\beta -2}\nonumber\\
 &:=& T^{4\alpha } S^{4\beta } U_{T,S}.\label{U}
\end{eqnarray}
Using the asymptotic behavior of the Bessel function (\ref{asymbessel}), we have that
\begin{equation*}
\frac{J_{0} \left(2 \sqrt{ \theta (t-u)(s-v)TS} \right)}{(TS) ^{-\frac{1}{4}+ \eps}}\to _{T,S\to \infty} 0
\end{equation*}
for almost every $t,u,s,v \in (0,1)$ and for every $\eps >0$. We will next apply the dominated convergence theorem. To this end, using again relation (\ref{asymbessel}), it suffices to show that the integral
\begin{eqnarray}
 I&=&\int_{0} ^{1}dt \int_{0}^{t} du\int_{0} ^{1}dt_{0} \int_{0}^{t_{0}} du_{0} \int_{0}^{1} ds\int_{0} ^{s}dv \int_{0}^{1} ds_{0} \int_{0}^{s_{0}} dv_{0} \nonumber \\
&& \left( (t-u)(s-v)(t_{0}-u_{0})(s_{0}-v_{0})\right) ^{-\frac{1}{4}} \nonumber \\
&& \vert t-t_{0} \vert ^{2\alpha -2} \vert u-u_{0} \vert ^{2\alpha -2} \vert s-s_{0} \vert ^{2\beta -2} \vert v-v_{0} \vert ^{2\beta -2}  \label{I}
\end{eqnarray}
is finite. This is proved in the following lemma.

 \qed

 \begin{remark}
In the above statement we can replace the normalization $T ^{-2\alpha + \frac{1}{4}-\eps } S^{-2\beta + \frac{1}{4}-\eps}$ by $T ^{-2\alpha + \frac{1}{4} } S^{-2\beta + \frac{1}{4}}f(T,S)$ where $f(T,S)$ is a deterministic function which converges to zero as $T,S\to \infty$. This is a consequence of the proof below.
\end{remark}

\begin{lemma}\label{ln1}
Let $I$ be given by (\ref{I}). Then for $\frac{1}{2}<\alpha, \beta <\frac{5}{8}$ the integral $I$ is finite.
\end{lemma}
{\bf  Proof: } Consider the integral
\begin{eqnarray*}
&& \int_{0}^{s_{0}} dv_{0} \left( (t-u)(s-v)\right) ^{-\frac{1}{4}} \left( (t_{0}-u_{0})(s_{0}-v_{0})\right) ^{-\frac{1}{4}} \\
&&	\times \vert t-t_{0} \vert ^{2\alpha -2} \vert u-u_{0} \vert ^{2\alpha -2} \vert s-s_{0} \vert ^{2\beta -2} \vert v-v_{0} \vert ^{2\beta -2}  \nonumber \\
&=& \left( (t-u)(s-v)(t_{0}-u_{0})\right) ^{-\frac{1}{4}} \vert t-t_{0} \vert ^{2\alpha -2} \vert u-u_{0} \vert ^{2\alpha -2} \vert s-s_{0} \vert ^{2\beta -2} \\
&& \times \int_{0}^{s_{0}} dv_{0} \vert v-v_{0} \vert ^{2\beta -2} \left(s_{0}-v_{0}\right) ^{-\frac{1}{4}}.
\end{eqnarray*}
If $s_{0} < v$ we have
\begin{equation*}
\int_{0}^{s_{0}} dv_{0} \vert v-v_{0} \vert ^{2\beta -2} \left(s_{0}-v_{0}\right) ^{-\frac{1}{4}} = \int_{0}^{s_{0}} dv_{0} \left( v-v_{0} \right)^{2\beta -2} \left(s_{0}-v_{0}\right) ^{-\frac{1}{4}}
\end{equation*}
making the change of variables $z=\frac{s_{0}-v_{0}}{v-v_{0}}$ we get
\begin{eqnarray*}
&&\int_{0}^{\frac{s_{0}}{v}} \left( \frac{v-s_{0}}{1-z} \right)^{2\beta -2} \left( \frac{z(v-s_{0})}{1-z} \right)^{-\frac{1}{4}} \left( \frac{v-s_{0}}{(1-z)^{2}} \right) dz \\
&& \leq \left(v-s_{0}\right)^{2\beta - \frac{5}{4} } \int_{0}^{1} z^{-\frac{1}{4}} (1-z)^{\frac{1}{4}-2\beta} dz = \left( v-s_{0} \right)^{2\beta - \frac{5}{4}} \tilde{\beta} \left( \frac{3}{4}, \frac{5}{4}-2\beta \right)
\end{eqnarray*}
where $\tilde{\beta}$ is the \emph{Beta} function. The expresion above is finite for $\beta<\frac{5}{8}$.\\
Now, if $v\leq s_{0}$
\begin{eqnarray*}
&&\int_{0}^{s_{0}} dv_{0} \vert v-v_{0} \vert ^{2\beta -2} \left(s_{0}-v_{0}\right) ^{-\frac{1}{4}} \\
&&= \int_{0}^{v} dv_{0} \left( v-v_{0} \right)^{2\beta -2} \left(s_{0}-v_{0}\right) ^{-\frac{1}{4}} + \int_{v}^{s_{0}} dv_{0} \left(v_{0}-v \right)^{2\beta -2} \left(s_{0}-v_{0}\right) ^{-\frac{1}{4}}.
\end{eqnarray*}
For the first integral in the right hand side we make the change of variables $z=\frac{v-v_{0}}{s_{0}-v_{0}}$, and for the second one we make $z=\frac{s_{0}-v_{0}}{v_{0}-v}$ Then we get
\begin{eqnarray*}
&\leq & \left( s_{0} -v \right)^{2\beta - \frac{5}{4}} \left( \int_{0}^{1} z^{2\beta-2} (1-z)^{\frac{1}{4}-2\beta} dz + \int_{0}^{\infty} z^{-\frac{1}{4}} (1+z)^{\frac{1}{4}-2\beta} dz \right) \nonumber \\
&=& \left( s_{0} -v \right)^{2\beta - \frac{5}{4}} \left( \tilde{\beta} \left( 2\beta -1, \frac{5}{4}-2\beta \right) + \frac{\Gamma(\frac{3}{4})\Gamma(2\beta -1)}{\Gamma(2\beta -\frac{1}{4})}\ _{2}F_{1}(0,3/4,2\beta -1/4;0) \right) \nonumber \\
&=& \left( s_{0} -v \right)^{2\beta - \frac{5}{4}} \left( \tilde{\beta} \left( 2\beta -1, \frac{5}{4}-2\beta \right) + \tilde{\beta} \left( \frac{3}{4}, 2\beta - 1 \right)  \ _{2}F_{1}(0,3/4,2\beta -1/4;0) \right) 
\end{eqnarray*}
wich is finite for $\frac{1}{2}<\alpha, \beta <\frac{5}{8}$. Here $\Gamma$ is the \emph{Gamma} function, $ _{2}F_{1}$ is the \emph{Hypergeometric} function, and we have make use of the property $\tilde{\beta}(x,y)=\frac{\Gamma(x) \Gamma(y)}{\Gamma(x+y)} $ and the formula (see \cite{S}, formula 1.6.7)
\begin{equation*}
_{2}F_{1}\left( a,b,c;1-x \right) = \frac{\Gamma(c)}{\Gamma(b)\Gamma(c-b)} \int_{0}^{\infty} w^{b-1}(1+w)^{a-c}(1+wx)^{-a} dw.
\end{equation*}


Proceeding in a similar way with the other seven integrals we conclude that $I$ is finite.

 \qed

\begin{remark}
We may note that the prescence of the $cos$ function in the asymptotic behavior of the Bessel function does not allow to obtain a renormalization for $F_{T,S}$ in terms of the powers of $T$ and $S$ like in \cite{HN}.
\end{remark}

The following proposition is a consequence of the Proposition \ref{pn1} and of the hypercontractivity property of multiple stochastic integrals.
\begin{prop}\label{pn2}
For every $\eps >0$ and for $\frac{1}{2}<\alpha, \beta <\frac{5}{8}$ the sequence
$$\frac{1}{ \sqrt{ T^{4\alpha -\frac{1}{2}+ \eps } S^{4\beta -\frac{1}{2} + \eps }} } \int_{0}^{T} \int_{0}^{S} X_{t,s} dB^{\alpha, \beta}_{t,s}$$
converges to zero a.s. as $T,S\to \infty$.
\end{prop}
{\bf Proof: } As in \cite{HN} we can replace the couple $(T,S)$ by  a discrete sequence  $(T_{M}, S_{N})$ such that $T_{M}, S_{N}$ converge to infinity as $M, N\to \infty.$ This is possible since the nominator and denominator in (\ref{theta-theta}) are continuous a.s. with respect to $(T,S)$. Indeed, the fact that the divergence integral in the nominator is continuous follows from \cite{N}, page 293 since the integrand $X$ is regular enough, and the integral $dtds$ in the nominator is clearly continuous a.s. with respect to the couple $(T,S)$. For simplicity, we will assume that $(T_{M}, S_{N})=(M,N)$.  Let us show that
$$A_{M,N} := M^{-2\alpha + \frac{1}{4} -\eps} N^{-2\beta + \frac{1}{4} -\eps} \int_{0}^{M} \int_{0}^{N}  X_{t,s} dB^{\alpha, \beta}_{t,s}$$
converges to zero a.s. as $M,N$ tend to infinity. We will use the Borel-Cantelli lemma. To do this, we will estimate $P(A_{M,N} >(MN) ^{-\gamma} )$ for some $\gamma >0$. For every $p\geq 1$ we have
\begin{equation*}
P(A_{M,N} >(MN) ^{-\gamma} )\leq \left( MN \right) ^{p\gamma} \mathbf{E}| A_{M,N} | ^{p}
\end{equation*}
and since $A_{M,N}$ is a multiple integral in the second Wiener chaos, the inequality (\ref{hyper}) and Proposition \ref{pn1} implies that
\begin{equation*}
\mathbf{E}| A_{M,N} | ^{p} \leq c(p) \left( \mathbf{E} A_{M,N} ^{2} \right) ^{\frac{p}{2} }\leq c(p, \alpha, \beta) (MN) ^{-\eps p}.
\end{equation*}
Putting together the two above bounds, we get
$$\sum_{M> M_{0}, N>N_{0} } P(A_{M,N} >(MN) ^{-\gamma} )\leq c(p, \alpha, \beta) \sum_{M> M_{0}, N>N_{0} } (MN) ^{p(\gamma -\eps)} $$
and this series is convergent when
$$(\eps-\gamma)p>1 \ \ \mbox{or equivalently } \ \ \gamma < \eps -\frac{1}{p}.$$
For every given $\eps >0$ and for $p$ large enough we can always chose a real number $\gamma $ such that $0<\gamma < \eps -\frac{1}{p}$. This, together with the Borel-Cantelli lemma allows us to finish the proof. \qed

\vskip0.3cm
The next step is to analyze the denominator in formula (\ref{theta-theta}). We have the following estimate.

\begin{prop}\label{pn3}
For any $\eps >0$ and for any $\alpha, \beta \in (\frac{1}{2}, \frac{5}{8})$
\begin{equation*}
\frac{1}{ T ^{2\alpha + \frac{1}{2} -\eps} S^{2\beta + \frac{1}{2} -\eps}}\mathbf{E} \int_{0} ^{T}dt \int_{0} ^{S} ds X_{t,s } ^{2}\to _{T,S\to \infty} \infty.
\end{equation*}
\end{prop}
{\bf Proof: }By the isometry of multiple integrals (\ref{iso}), the equation (\ref{pe}) and the expression of the solution to (\ref{1})
\begin{eqnarray*}
\mathbf{E}X_{t,s} ^{2}&=& c(\alpha, \beta) \int_{0} ^{t} du \int_{0} ^{t} du_{0} \int_{0} ^{s}dv \int_{0} ^{s} dv_{0} \\
&&\times J_{0} ( 2\sqrt{ \theta (t-u) (s-v) })
J_{0} ( 2\sqrt{ \theta (t-u_{0}) (s-v_{0}) }) \vert u-u_{0} \vert ^{2\alpha -2} \vert v-v_{0} \vert ^{2\beta -2},
\end{eqnarray*}
making the change of variables $\tilde{t}= \frac{t}{T}, \tilde{s}=\frac{s}{S}$ and similar for the other variables we can write
\begin{eqnarray*}
\int_{0}^{T} \int_{0} ^{S} dsdt \mathbf{E}X_{t,s} ^{2}&=& c(\alpha, \beta)T^{2\alpha +1} S^{2\beta +1}\int_{0}^{1} dt \int_{0} ^{1} ds
\int_{0} ^{t} du \int_{0} ^{t} du_{0} \int_{0} ^{s}dv \int_{0} ^{s} dv_{0} \\
&& \times J_{0} (2 \sqrt{ \theta (t-u) (s-v) TS})
J_{0} (2 \sqrt{ \theta (t-u_{0}) (s-v_{0})TS }) \\
&& \times \vert u-u_{0} \vert ^{2\alpha -2} \vert v-v_{0} \vert ^{2\beta -2}
\end{eqnarray*}
and thus
\begin{eqnarray*}
\frac{\int_{0}^{T} \int_{0} ^{S} dsdt \mathbf{E}X_{t,s} ^{2}}{ T ^{2\alpha + \frac{1}{2} -\eps} S^{2\beta + \frac{1}{2} -\eps}}
&=& c(\alpha, \beta) \int_{0}^{1} dt  \int_{0} ^{1} ds
\int_{0} ^{t} du \int_{0} ^{t} du_{0} \int_{0} ^{s}dv \int_{0} ^{s} dv_{0} \\
&& \times \frac{J_{0} (2 \sqrt{ \theta (t-u) (s-v) TS})}{(TS)^{-\frac{1}{4}-\frac{\eps}{2}}}
\frac{J_{0} (2 \sqrt{ \theta (t-u_{0}) (s-v_{0})TS })}{(TS) ^{-\frac{1}{4}-\frac{\eps}{2}}} \\
&& \times \vert u-u_{0} \vert ^{2\alpha -2} \vert v-v_{0} \vert ^{2\beta -2}.
\end{eqnarray*}
We will use the same idea as in the proof of Proposition \ref{pn1}. We notice first that
$$\frac{ J_{0} (2 \sqrt{ \theta (t-u) (s-v) TS})}{(TS)^{-\frac{1}{4} - \frac{\eps}{2}}}$$
converges to infinity as $T,S \to \infty$ by using (\ref{asymbessel}); also for $T,S $ large enough
$$\frac{ J_{0} (2 \sqrt{ \theta (t-u) (s-v) TS})}{(TS)^{-\frac{1}{4} }}$$
is  bounded by $c  \left( \vert t-u\vert \vert s-v\vert \right) ^{-\frac{1}{4}}$ almost everywhere $s,t,u,v$. Also we note that  the integral
\begin{eqnarray*}
&&\int_{0}^{1} dt  \int_{0} ^{1} ds
\int_{0} ^{t} du \int_{0} ^{t} du_{0} \int_{0} ^{s}dv \int_{0} ^{s} dv_{0}
 \vert u-u_{0} \vert ^{2\alpha -2} \vert v-v_{0} \vert ^{2\beta -2}\\
 &&\left( \vert t-u\vert \vert s-v\vert \right) ^{-\frac{1}{4}}\left( \vert t-u_{0}\vert \vert s-v_{0}\vert \right) ^{-\frac{1}{4}}
 \end{eqnarray*}
 is finite for $\alpha , \beta \in (\frac{1}{2}, \frac{5}{8})$ by using the same computations as in the proof of Lemma \ref{ln1}. \\
Using Fatou's lemma  (we use the following version of the Fatou's lemma: if $f_{n}$ is a sequence of functions such that $f_{n} \geq -g$ where $g$ is positive and integrable, then $\liminf \int f_{n} \geq \int \liminf f_{n}$) this implies that
\begin{eqnarray*}
\liminf_{T,S\to \infty} \frac{\int_{0}^{T} \int_{0} ^{S} dsdt \mathbf{E}X_{t,s} ^{2}}{ T ^{2\alpha + \frac{1}{2} -\eps} S^{2\beta + \frac{1}{2} -\eps}}
&\geq & c(\alpha,\beta)  \int_{0}^{1} dt  \int_{0} ^{1} ds
\int_{0} ^{t} du \int_{0} ^{t} du_{0} \int_{0} ^{s}dv \int_{0} ^{s} dv_{0}  \\
&\times & \liminf_{T,S\to \infty} \frac{J_{0} (2 \sqrt{ \theta (t-u) (s-v) TS})}{(TS)^{-\frac{1}{4}-\frac{\eps}{2}}}
\frac{J_{0} (2 \sqrt{ \theta (t-u_{0}) (s-v_{0})TS })}{(TS) ^{-\frac{1}{4}-\frac{\eps}{2}}} \\
&\times & \vert u-u_{0} \vert ^{2\alpha -2} \vert v-v_{0} \vert ^{2\beta -2} \ \ =\ \ \infty .
\end{eqnarray*}\qed

\vskip0.3cm

At this point we will need the following auxiliary lemma.
\begin{lemma}\label{ln2}
Consider a sequence of random variables $(A_{N})_{N}$ such that $\sum_{N} P(A_{N}>cN^{-\gamma} ) <\infty$ for some $\gamma >0$ and for every $c>0$ (which implies  $A_{N} \to 0 $ almost surely as $N\to \infty$.) Also consider a sequence of a.s. strictly positive random variables $(B_{N})_{N}$ such that $\mathbf{E}B_{N}\to \infty$ as $N\to \infty$. Then
\begin{equation*}
\frac{A_{N} }{B_{N}} \to_{N\to \infty} 0  \mbox{ almost surely}
\end{equation*}
\end{lemma}
{\bf Proof: }We will use again the Borel-Cantelli lemma. Let $C >0$ be arbitrary. Then
\begin{eqnarray*}
\sum_{N} P \left( \frac{A_{N}}{B_{N}}>N^{-\gamma } \right) &=&\sum_{N} P \left( \frac{A_{N}}{B_{N}}>N^{-\gamma }, B_{N}>C \right)+ \sum_{N} P\left(\frac{A_{N}}{B_{N}}>N^{-\gamma }, B_{N}<C\right)\\
&\leq & \sum_{N} P\left(\frac{A_{N}}{C}>N^{-\gamma }, B_{N}>C \right) + \sum_{N} P\left(\frac{A_{N}}{B_{N}}>N^{-\gamma }, B_{N}<C\right).
\end{eqnarray*}
Using the fact that $\mathbf{E}B_{N}\to \infty$ as $N\to \infty$ we obtain that
$$ P\left(\frac{A_{N}}{B_{N}}>N^{-\gamma }, B_{N}<C\right)=0$$
for $N $ large enough. By assumption, $ \sum_{N} P\left(\frac{A_{N}}{C}>N^{-\gamma }, B_{N}>C \right)<\infty$. Therefore
$$\sum_{N} P\left(\frac{A_{N}}{B_{N}}>N^{-\gamma }\right) <\infty$$ and therefore the conclusion follows.\\

 \qed

\begin{remark}
Lemma \ref{ln2} can be extended without difficulty to two-parameter sequences. That is, if $A_{M,N}$ is a sequence of random variables such that $\sum_{M,N\geq 0} P(A_{M,N}>c(MN)^{-\gamma} ) <\infty$ and $B_{M,N}$ is a sequence of positive random variables such that $\mathbf{E}B_{M,N} \to _{M,N\to \infty} \infty$ then
$$\frac{A_{M,N}}{B_{M,N}} \to _{M,N\to \infty} 0 \mbox{ almost surely}. $$
\end{remark}
\vskip0.3cm

 Let us state the main result of this section.
\begin{theorem}Let $\theta _{T,S}$ be the estimator given by (\ref{est1}). Suppose that $\alpha , \beta \in (\frac{1}{2}, \frac{5}{8})$. Then $\theta _{T,S}$ is a strongly consistent estimator for the parameter $\theta$, that is,
\begin{equation*}
\widehat{\theta} _{T,S} \to \theta \mbox{ almost surely as } T,S\to \infty.
\end{equation*}
\end{theorem}
{\bf Proof: }From relation (\ref{theta-theta}), the difference between the estimator and the true parameter is
\begin{eqnarray*}
\widehat{\theta}_{T}-\theta &=&\frac{ T^{-2\alpha + \frac{1}{4} -\eps } S ^{ -2\beta + \frac{1}{4} -\eps} \int_{0} ^{T} \int_{0}^{S} X_{t,s} dB^{\alpha , \beta } _{t,s}} { T^{-2\alpha + \frac{1}{4} -\eps } S ^{ -2\beta + \frac{1}{4} -\eps}\int_{0}^{T} dt \int_{0} ^{S} ds X_{t,s} ^{2}}\\
&=& \frac{ T^{-2\alpha + \frac{1}{4} -\eps } S ^{ -2\beta + \frac{1}{4} -\eps} \int_{0} ^{T} \int_{0}^{S} X_{t,s} dB^{\alpha , \beta } _{t,s}} {T ^{\frac{3}{4} -2\eps} S ^{\frac{3}{4} -2\eps} T^{-2\alpha - \frac{1}{2} +\eps} S^{-2\beta -\frac{1}{2} +\eps} \int_{0}^{T} dt \int_{0} ^{S} ds X_{t,s} ^{2}}.
\end{eqnarray*}The result is obtained by using Proposition \ref{pn2}, Proposition \ref{pn3} and Lemma \ref{ln2} (and the remark that follows after this lemma). \qed

\section{Asymptotic non-normality of the estimator}

We proved in the previous section that the least square estimator $\widehat{\theta}_{T,S}$ is strongly consistent. This property has been proved in the one-dimensional case in (\cite{HN}). Nevertheless, we show in this paragraph that the limiting distribution of the estimator is not the same in the one-parameter and two-parameter cases. This different behavior is somehow expected since the fractional Ornstein-Uhlenbeck sheet does not keep the properties of the fractional Ornstein-Uhlenbeck process (for example, the kernel $f$ given by (\ref{f}) can take any real value while in the one parameter case it is positive since it is given by an exponential function). In order to notice the difference between the one-parameter and the two-parameter case, we will focus only on the nominator in the right hand side of (\ref{theta-theta}).  To check the asymptotic normality  we use the following criterium:
(see Theorem 4 in \cite{NOT}).

\begin{theorem}
Let $(F_{k}, k\geq 1)$, $F_{k}=I_{n}(f_{k})$ (with $f_{k}\in
{\cal{H}}^{\odot n}$  for every $k\ge 1$) be a sequence of square
integrable random variables in the $n$ th Wiener chaos such that
$\mathbf{E}[F_{k}^{2}]\rightarrow 1$ as $k\rightarrow\infty.$ Then
the following are equivalent:

\begin{description}
\item {i)} The sequence $(F_{k}) _{k\geq0} $ converges in distribution to the normal law
${\cal {N}} (0,1)$.

\item {ii) }$\Vert DF_{k}\Vert_{{\mathcal{H}}}^{2}$ converges to $n$ in
$L^{2}(\Omega)$ as $k\rightarrow\infty$.
\end{description}
\end{theorem}

\vskip0.3cm

Denote by
\begin{equation*}
\sigma ^{2} _{T,S}= \mathbf{E} \left(  \int_{0} ^{T} \int_{0}^{S} X_{t,s} dB^{\alpha, \beta } _{t,s}\right) ^{2} = \mathbf{E} \left( F_{T,S} \right)^{2}.
\end{equation*}

\vskip0.5cm
\begin{prop}
Assume $\frac{1}{2} <\alpha , \beta <\frac{5}{8}$  and let $F_{T,S}$ be given  by (\ref{FTS}). Then when $T,S$ tend to infinity,
$$\frac{1}{\sigma _{T,S}}F_{T,S}$$
 does not converges in distribution to the normal law $N(0, 1)$.

\end{prop}
{\bf Proof: }We will prove that $\Vert D\frac{1}{\sigma _{T,S}}F_{T,S} \Vert _{{\cal{H}}^{\alpha, \beta }} ^{2}$
does not converges to $2$ in $L^{2}(\Omega)$ as $T,S \to \infty$.

Let us denote by $\tilde{f}$ the symmetrization of $f$ with respect
to the variables $(t,s), (u,v)$
\begin{equation*}
\tilde{f} \left( (t,s), (u,v)\right) =\frac{1}{2} \left( f(t,s,u,v)+
f(u,v,t,s) \right).
\end{equation*}

We have
\begin{eqnarray*}
 D_{t,s} \frac{1}{\sigma _{T,S}} F_{T,S}&=&\frac{1}{\sigma _{T,S}}  D_{t,s} I_{2} \left(f \right) = \frac{1}{\sigma _{T,S}}  D_{t,s} I_{2} \left( \tilde{f}
 \right) =\frac{1}{\sigma _{T,S}}\left( I_{1}( f(\cdot, \cdot, t,s)) +
 I_{1}(t,s, \cdot, \cdot ) \right).
\end{eqnarray*}
Here $"\cdot, \cdot "$ represents the variable with respect to which
the integral $I_{1}$ acts.

We obtain
\begin{eqnarray*}
\Vert D\frac{1}{\sigma _{T,S}}F_{T,S} \Vert _{{\cal{H}}^{\alpha, \beta }} ^{2}&=&
\frac{1}{\sigma^{2} _{T,S}} c(\alpha) c(\beta) \int_{0}^{T} dt \int_{0}^{S}ds
\int_{0}^{T} dt_{0} \int_{0} ^{S}ds_{0} \\
&&\times \left[ I_{1} \left( f(\cdot, \cdot, t,s)\right) + I_{1}
\left( f(t,s,\cdot, \cdot ) \right) \right]
 \left[ I_{1} \left( f(\cdot, \cdot, t_{0},s_{0})\right) + I_{1}
\left( f(t_{0},s_{0},\cdot, \cdot ) \right) \right]\\
&&\times  \vert t-t_{0}
\vert ^{2\alpha -2} \vert s-s_{0}\vert ^{2\beta -2} \\
&=&\frac{1}{\sigma ^{2} _{T,S}} c(\alpha) c(\beta) \int_{0}^{T} dt \int_{0}^{S}ds
\int_{0}^{T} dt_{0} \int_{0} ^{S}ds_{0} \\
&&\times \left[ I_{2} \left( f(\cdot, \cdot, t,s) \otimes f(\cdot,
\cdot, t_{0}, s_{0}) \right) +I_{2} \left( f(\cdot, \cdot, t,s)
\otimes f(t_{0}, s_{0}, \cdot, \cdot) \right)\right. \\
&&\left. +  I_{2}\left( f(t,s, \cdot, \cdot) \otimes f(\cdot, \cdot,
t_{0}, s_{0}) \right)+I_{2}\left( f(t,s, \cdot, \cdot) \otimes f(
t_{0}, s_{0}, \cdot, \cdot) \right)\right]\\
&&\times  \vert t-t_{0}
\vert ^{2\alpha -2} \vert s-s_{0} \vert ^{2\beta -2} +\mathbf{E}\Vert D\frac{1}{\sigma _{T,S}}F_{T,S} \Vert _{{\cal{H}}^{\alpha, \beta }} ^{2}\\
&:=& A^{(1)}_{T,S} +A^{(2)}_{T,S} +A^{(3)}_{T,S}
+A^{(4)}_{T,S}+\mathbf{E}\Vert D\frac{1}{\sigma _{T,S}}F_{T,S}  \Vert _{{\cal{H}}^{\alpha, \beta }}
^{2}
\end{eqnarray*}
Let us note that that $\mathbf{E}\Vert D\frac{1}{\sigma _{T,S}}F_{T,S} \Vert _{{\cal{H}}^{\alpha,
\beta }} ^{2}$ is equal to $2$. This follows from the fact
that for any multiple integral of order $n$ we have
$$\mathbf{E}\Vert DI_{n}(f)\Vert ^{2} _{{\cal{H}}} = n\mathbf{E}I_{n} (f) ^{2}.$$

It remains to show that the terms containing multiple integrals of
order 2 does not converges to zero in $L^{2}(\Omega)$ as $T,S\to \infty$. The
first and the fourth summand are similar, as they are the second one
and the third one. We will handle only the first summand denoted by
$A^{(1)}_{T,S}$ (because the other three terms can be studied analogously). We can write, using the definition of the scalar
product in the Hilbert space $\left( {\cal{H}}^{\alpha , \beta}
\right) ^{\otimes 2}$ and the expression (\ref{fjo}) of the kernel
$f$ in terms of the Bessel function $J_{0}$\\

\begin{eqnarray*}
\mathbf{E}\left| A^{(1)}_{T,S}\right| ^{2}&=&\frac{1}{\sigma ^{4}_{T,S}}\mathbf{E} \left( \int_{0}^{T} dt\int_{0}^{S} ds \int_{0}
^{T}dt_{0} \int_{0}^{S} ds_{0}\vert t-t_{0}\vert ^{2\alpha -2}\vert
s-s_{0}\vert ^{2\beta -2} \right. \\
&&\left. I_{2} \left( f(\cdot, \cdot, t,s) \otimes f(\cdot, \cdot,
t_{0}, s_{0}) \right) \ \ \right) ^{2}\\
&\sim &2\int_{0}^{T} dt\int_{0}^{S} ds \int_{0} ^{T}dt_{0}
\int_{0}^{S} ds_{0}\int_{0}^{T} du\int_{0}^{S} dv\int_{0} ^{T}du_{0}
\int_{0}^{S} dv_{0}\\
&&\vert t-t_{0}\vert ^{2\alpha -2}\vert s-s_{0}\vert ^{2\beta -2}
\vert u-u_{0}\vert ^{2\alpha -2}\vert v-v_{0}\vert ^{2\beta -2}\\
&&\langle f(\cdot, \cdot, t,s) \otimes f(\cdot, \cdot, t_{0}, s_{0})
, f(\cdot, \cdot, u,v) \otimes f(\cdot, \cdot, u_{0}, v_{0})\rangle
_{\left( {\cal{H}}^{\alpha , \beta} \right) ^{\otimes 2}}\\
&=&C\frac{1}{\sigma ^{4}_{T,S}}\int_{0}^{T} dt\int_{0}^{S} ds \int_{0} ^{T}dt_{0} \int_{0}^{S}
ds_{0}\int_{0}^{T} du\int_{0}^{S} dv\int_{0} ^{T}du_{0}
\int_{0}^{S} dv_{0}\\
&&\int_{0}^{T} dx\int_{0}^{S} dy \int_{0} ^{T}dx_{0} \int_{0}^{S}
dy_{0}\int_{0}^{T} da\int_{0}^{S} db\int_{0} ^{T}da_{0}
\int_{0}^{S} db_{0}\\
&&f(x,y,t,s) f(x_{0}, y_{0}, t_{0}, s_{0}) f(a,b,u,v) f(a_{0},
b_{0}, t_{0}, s_{0}) \\
&&\vert t-t_{0}\vert ^{2\alpha -2}\vert s-s_{0}\vert ^{2\beta -2}
\vert u-u_{0}\vert ^{2\alpha -2}\vert v-v_{0}\vert ^{2\beta -2}\\
&&\vert x-a\vert ^{2\alpha -2}\vert y-b\vert ^{2\beta -2}
\vert x_{0}-a_{0}\vert ^{2\alpha -2}\vert y_{0}-b_{0}\vert ^{2\beta -2}\\
&=&C\frac{1}{\sigma ^{4}_{T,S}}\int_{0}^{T} dx\int_{0}^{S} dy \int_{0} ^{T}dx_{0} \int_{0}^{S}
dy_{0}\int_{0}^{T} da\int_{0}^{S} db\int_{0} ^{T}da_{0}
\int_{0}^{S} db_{0}\\
&&\int_{0}^{x} dt\int_{0}^{y} ds \int_{0} ^{x_{0}}dt_{0}
\int_{0}^{y_{0}} ds_{0}\int_{0}^{a} du\int_{0}^{b} dv\int_{0}
^{a_{0}}du_{0}
\int_{0}^{b_{0}} dv_{0}\\
&& J_{0} (2\sqrt{\theta (x-t) (y-s) }) J_{0} (2\sqrt{\theta
(x_{0}-t_{0}) (y_{0}-s_{0}) })\\
 && J_{0} (2\sqrt{\theta (a-u) (b-v)
})J_{0} (2\sqrt{\theta
(a_{0}-u_{0}) (b_{0}-v_{0}) }) \\
 &&\vert t-t_{0}\vert ^{2\alpha
-2}\vert s-s_{0}\vert ^{2\beta -2}
\vert u-u_{0}\vert ^{2\alpha -2}\vert v-v_{0}\vert ^{2\beta -2}\\
&&\vert x-a\vert ^{2\alpha -2}\vert y-b\vert ^{2\beta -2}
\vert x_{0}-a_{0}\vert ^{2\alpha -2}\vert y_{0}-b_{0}\vert ^{2\beta -2}
\end{eqnarray*}
making $\tilde{t}= \frac{t}{T}, \tilde{s}= \frac{s}{S}$ and similar for the other variables we get
\begin{eqnarray*}
\mathbf{E}\left| A^{(1)}_{T,S}\right| ^{2}&=&
C\frac{1}{\sigma ^{4}_{T,S}}T^{8\alpha} S^{8\beta} \int_{0}^{1} dx\int_{0}^{1} dy \int_{0} ^{1}dx_{0} \int_{0}^{1}
dy_{0}\int_{0}^{1} da\int_{0}^{1} db\int_{0} ^{1}da_{0}
\int_{0}^{1} db_{0}\\
&&\int_{0}^{x} dt\int_{0}^{y} ds \int_{0} ^{x_{0}}dt_{0}
\int_{0}^{y_{0}} ds_{0}\int_{0}^{a} du\int_{0}^{b} dv\int_{0}
^{a_{0}}du_{0}
\int_{0}^{b_{0}} dv_{0}\\
&& J_{0} (2\sqrt{\theta (x-t) (y-s) TS}) J_{0} (2\sqrt{\theta
(x_{0}-t_{0}) (y_{0}-s_{0})TS })\\
 && J_{0} (2\sqrt{\theta (a-u) (b-v)
TS})J_{0} (2\sqrt{\theta
(a_{0}-u_{0}) (b_{0}-v_{0})TS }) \\
 &&\vert t-t_{0}\vert ^{2\alpha
-2}\vert s-s_{0}\vert ^{2\beta -2}
\vert u-u_{0}\vert ^{2\alpha -2}\vert v-v_{0}\vert ^{2\beta -2}\\
&&\vert x-a\vert ^{2\alpha -2}\vert y-b\vert ^{2\beta -2}
\vert x_{0}-a_{0}\vert ^{2\alpha -2}\vert y_{0}-b_{0}\vert ^{2\beta -2}\\.
\end{eqnarray*}
Using the asymptotic behavior of the Bessel function when its
variable is close to infinity (see (\ref{asymbessel})) we see that
\begin{eqnarray*}
\mathbf{E}\left| A^{(1)}_{T,S}\right| ^{2}&\approx &
C\frac{1}{\sigma ^{4}_{T,S}}T^{8\alpha} S^{8\beta} \int_{0}^{1} dx\int_{0}^{1} dy \int_{0} ^{1}dx_{0} \int_{0}^{1}
dy_{0}\int_{0}^{1} da\int_{0}^{1} db\int_{0} ^{1}da_{0}
\int_{0}^{1} db_{0}\\
&&\int_{0}^{x} dt\int_{0}^{y} ds \int_{0} ^{x_{0}}dt_{0}
\int_{0}^{y_{0}} ds_{0}\int_{0}^{a} du\int_{0}^{b} dv\int_{0}
^{a_{0}}du_{0}
\int_{0}^{b_{0}} dv_{0}\\
&&\left( (x-t)(y-s)(x_{0}-t_{0})(y_{0}-s_{0})(a-u) (b-v)(a_{0}-u_{0}) (b_{0}-v_{0}) \right)^{\frac{-1}{4}}\\
&&\vert t-t_{0}\vert ^{2\alpha
-2}\vert s-s_{0}\vert ^{2\beta -2}
\vert u-u_{0}\vert ^{2\alpha -2}\vert v-v_{0}\vert ^{2\beta -2}\\
&&\vert x-a\vert ^{2\alpha -2}\vert y-b\vert ^{2\beta -2}
\vert x_{0}-a_{0}\vert ^{2\alpha -2}\vert y_{0}-b_{0}\vert ^{2\beta -2}.
\end{eqnarray*}
and considering the fact that the last integral
is finite for $\frac{1}{2}<\alpha , \beta < \frac{5}{8}$ (the proof of this fact is similar to the proof of Lemma \ref{ln1} ) it is straightforward to see that, for $T,S$ close to infinity, the quantity $\mathbf{E}\left| A^{(1)}_{T,S}\right| ^{2}$ does not converges to zero. \\

 \qed

{\bf Acknowledgement: } The authors would like to thank Prof. Soledad Torres for interesting discussions.

\end{document}